\catcode`\@=11

\font\teneufm=eufm10
\font\seveneufm=eufm7
\font\fiveeufm=eufm5
\newfam\eufmfam
\textfont\eufmfam=\teneufm  \scriptfont\eufmfam=\seveneufm
  \scriptscriptfont\eufmfam=\fiveeufm
\def\eufm@{\hexnumber@\eufmfam}

\font\tenmsa=msam10
\font\sevenmsa=msam7
\font\fivemsa=msam5
\font\tenmsb=msbm10
\font\sevenmsb=msbm7
\font\fivemsb=msbm5
\newfam\msafam
\newfam\msbfam
\textfont\msafam=\tenmsa  \scriptfont\msafam=\sevenmsa
  \scriptscriptfont\msafam=\fivemsa
\textfont\msbfam=\tenmsb  \scriptfont\msbfam=\sevenmsb
  \scriptscriptfont\msbfam=\fivemsb

\def\hexnumber@#1{\ifnum#1<10 \number#1\else
 \ifnum#1=10 A\else\ifnum#1=11 B\else\ifnum#1=12 C\else
 \ifnum#1=13 D\else\ifnum#1=14 E\else\ifnum#1=15 F\fi\fi\fi\fi\fi\fi\fi}

\def\msa@{\hexnumber@\msafam}
\def\msb@{\hexnumber@\msbfam}

\mathchardef\gx="2\eufm@78
\mathchardef\gg="2\eufm@67
\mathchardef\gm="2\eufm@6D
\mathchardef\gn="2\eufm@6E
\mathchardef\gd="2\eufm@64

\mathchardef\boxdot="2\msa@00
\mathchardef\boxplus="2\msa@01
\mathchardef\boxtimes="2\msa@02
\mathchardef\square="0\msa@03
\mathchardef\blacksquare="0\msa@04
\mathchardef\centerdot="2\msa@05
\mathchardef\lozenge="0\msa@06
\mathchardef\blacklozenge="0\msa@07
\mathchardef\circlearrowright="3\msa@08
\mathchardef\circlearrowleft="3\msa@09
\mathchardef\rightleftharpoons="3\msa@0A
\mathchardef\leftrightharpoons="3\msa@0B
\mathchardef\boxminus="2\msa@0C
\mathchardef\Vdash="3\msa@0D
\mathchardef\Vvdash="3\msa@0E
\mathchardef\vDash="3\msa@0F
\mathchardef\twoheadrightarrow="3\msa@10
\mathchardef\twoheadleftarrow="3\msa@11
\mathchardef\leftleftarrows="3\msa@12
\mathchardef\rightrightarrows="3\msa@13
\mathchardef\upuparrows="3\msa@14
\mathchardef\downdownarrows="3\msa@15
\mathchardef\upharpoonright="3\msa@16

\mathchardef\downharpoonright="3\msa@17
\mathchardef\upharpoonleft="3\msa@18
\mathchardef\downharpoonleft="3\msa@19
\mathchardef\rightarrowtail="3\msa@1A
\mathchardef\leftarrowtail="3\msa@1B
\mathchardef\leftrightarrows="3\msa@1C
\mathchardef\rightleftarrows="3\msa@1D
\mathchardef\Lsh="3\msa@1E
\mathchardef\Rsh="3\msa@1F
\mathchardef\rightsquigarrow="3\msa@20
\mathchardef\leftrightsquigarrow="3\msa@21
\mathchardef\looparrowleft="3\msa@22
\mathchardef\looparrowright="3\msa@23
\mathchardef\circeq="3\msa@24
\mathchardef\succsim="3\msa@25
\mathchardef\gtrsim="3\msa@26
\mathchardef\gtrapprox="3\msa@27
\mathchardef\multimap="3\msa@28
\mathchardef\therefore="3\msa@29
\mathchardef\because="3\msa@2A
\mathchardef\doteqdot="3\msa@2B

\mathchardef\triangleq="3\msa@2C
\mathchardef\precsim="3\msa@2D
\mathchardef\lesssim="3\msa@2E
\mathchardef\lessapprox="3\msa@2F
\mathchardef\eqslantless="3\msa@30
\mathchardef\eqslantgtr="3\msa@31
\mathchardef\curlyeqprec="3\msa@32
\mathchardef\curlyeqsucc="3\msa@33
\mathchardef\preccurlyeq="3\msa@34
\mathchardef\leqq="3\msa@35
\mathchardef\leqslant="3\msa@36
\mathchardef\lessgtr="3\msa@37
\mathchardef\backprime="0\msa@38
\mathchardef\risingdotseq="3\msa@3A
\mathchardef\fallingdotseq="3\msa@3B
\mathchardef\succcurlyeq="3\msa@3C
\mathchardef\geqq="3\msa@3D
\mathchardef\geqslant="3\msa@3E
\mathchardef\gtrless="3\msa@3F
\mathchardef\sqsubset="3\msa@40
\mathchardef\sqsupset="3\msa@41
\mathchardef\trianglerighteq="3\msa@44
\mathchardef\trianglelefteq="3\msa@45
\mathchardef\bigstar="0\msa@46
\mathchardef\between="3\msa@47
\mathchardef\blacktriangledown="0\msa@48
\mathchardef\blacktriangleright="3\msa@49
\mathchardef\blacktriangleleft="3\msa@4A
\mathchardef\blacktriangle="0\msa@4E
\mathchardef\triangledown="0\msa@4F
\mathchardef\eqcirc="3\msa@50
\mathchardef\lesseqgtr="3\msa@51
\mathchardef\gtreqless="3\msa@52
\mathchardef\lesseqqgtr="3\msa@53
\mathchardef\gtreqqless="3\msa@54
\mathchardef\Rrightarrow="3\msa@56
\mathchardef\Lleftarrow="3\msa@57
\mathchardef\veebar="2\msa@59
\mathchardef\barwedge="2\msa@5A
\mathchardef\doublebarwedge="2\msa@5B
\mathchardef\angle="0\msa@5C
\mathchardef\measuredangle="0\msa@5D
\mathchardef\sphericalangle="0\msa@5E
\mathchardef\varpropto="3\msa@5F
\mathchardef\smallsmile="3\msa@60
\mathchardef\smallfrown="3\msa@61
\mathchardef\Subset="3\msa@62
\mathchardef\Supset="3\msa@63
\mathchardef\Cup="2\msa@64

\mathchardef\Cap="2\msa@65

\mathchardef\curlywedge="2\msa@66
\mathchardef\curlyvee="2\msa@67
\mathchardef\leftthreetimes="2\msa@68
\mathchardef\rightthreetimes="2\msa@69
\mathchardef\subseteqq="3\msa@6A
\mathchardef\supseteqq="3\msa@6B
\mathchardef\bumpeq="3\msa@6C
\mathchardef\Bumpeq="3\msa@6D
\mathchardef\lll="3\msa@6E

\mathchardef\ggg="3\msa@6F

\mathchardef\circledS="0\msa@73
\mathchardef\pitchfork="3\msa@74
\mathchardef\dotplus="2\msa@75
\mathchardef\backsim="3\msa@76
\mathchardef\backsimeq="3\msa@77
\mathchardef\complement="0\msa@7B
\mathchardef\intercal="2\msa@7C
\mathchardef\circledcirc="2\msa@7D
\mathchardef\circledast="2\msa@7E
\mathchardef\circleddash="2\msa@7F
\def\ulcorner{\delimiter"4\msa@70\msa@70 }
\def\urcorner{\delimiter"5\msa@71\msa@71 }
\def\llcorner{\delimiter"4\msa@78\msa@78 }
\def\lrcorner{\delimiter"5\msa@79\msa@79 }
\def\yen{\mathhexbox\msa@55 }
\def\checkmark{\mathhexbox\msa@58 }
\def\circledR{\mathhexbox\msa@72 }
\def\maltese{\mathhexbox\msa@7A }
\mathchardef\lvertneqq="3\msb@00
\mathchardef\gvertneqq="3\msb@01
\mathchardef\nleq="3\msb@02
\mathchardef\ngeq="3\msb@03
\mathchardef\nless="3\msb@04
\mathchardef\ngtr="3\msb@05
\mathchardef\nprec="3\msb@06
\mathchardef\nsucc="3\msb@07
\mathchardef\lneqq="3\msb@08
\mathchardef\gneqq="3\msb@09
\mathchardef\nleqslant="3\msb@0A
\mathchardef\ngeqslant="3\msb@0B
\mathchardef\lneq="3\msb@0C
\mathchardef\gneq="3\msb@0D
\mathchardef\npreceq="3\msb@0E
\mathchardef\nsucceq="3\msb@0F
\mathchardef\precnsim="3\msb@10
\mathchardef\succnsim="3\msb@11
\mathchardef\lnsim="3\msb@12
\mathchardef\gnsim="3\msb@13
\mathchardef\nleqq="3\msb@14
\mathchardef\ngeqq="3\msb@15
\mathchardef\precneqq="3\msb@16
\mathchardef\succneqq="3\msb@17
\mathchardef\precnapprox="3\msb@18
\mathchardef\succnapprox="3\msb@19
\mathchardef\lnapprox="3\msb@1A
\mathchardef\gnapprox="3\msb@1B
\mathchardef\nsim="3\msb@1C
\mathchardef\napprox="3\msb@1D
\mathchardef\nsubseteqq="3\msb@22
\mathchardef\nsupseteqq="3\msb@23
\mathchardef\subsetneqq="3\msb@24
\mathchardef\supsetneqq="3\msb@25
\mathchardef\subsetneq="3\msb@28
\mathchardef\supsetneq="3\msb@29
\mathchardef\nsubseteq="3\msb@2A
\mathchardef\nsupseteq="3\msb@2B
\mathchardef\nparallel="3\msb@2C
\mathchardef\nmid="3\msb@2D
\mathchardef\nshortmid="3\msb@2E
\mathchardef\nshortparallel="3\msb@2F
\mathchardef\nvdash="3\msb@30
\mathchardef\nVdash="3\msb@31
\mathchardef\nvDash="3\msb@32
\mathchardef\nVDash="3\msb@33
\mathchardef\ntrianglerighteq="3\msb@34
\mathchardef\ntrianglelefteq="3\msb@35
\mathchardef\ntriangleleft="3\msb@36
\mathchardef\ntriangleright="3\msb@37
\mathchardef\nleftarrow="3\msb@38
\mathchardef\nrightarrow="3\msb@39
\mathchardef\nLeftarrow="3\msb@3A
\mathchardef\nRightarrow="3\msb@3B
\mathchardef\nLeftrightarrow="3\msb@3C
\mathchardef\nleftrightarrow="3\msb@3D
\mathchardef\divideontimes="2\msb@3E
\mathchardef\varnothing="0\msb@3F
\mathchardef\nexists="0\msb@40
\mathchardef\mho="0\msb@66
\mathchardef\thorn="0\msb@67
\mathchardef\beth="0\msb@69
\mathchardef\gimel="0\msb@6A
\mathchardef\daleth="0\msb@6B
\mathchardef\lessdot="3\msb@6C
\mathchardef\gtrdot="3\msb@6D
\mathchardef\ltimes="2\msb@6E
\mathchardef\rtimes="2\msb@6F
\mathchardef\shortmid="3\msb@70
\mathchardef\shortparallel="3\msb@71
\mathchardef\smallsetminus="2\msb@72
\mathchardef\thicksim="3\msb@73
\mathchardef\thickapprox="3\msb@74
\mathchardef\approxeq="3\msb@75
\mathchardef\succapprox="3\msb@76
\mathchardef\precapprox="3\msb@77
\mathchardef\curvearrowleft="3\msb@78
\mathchardef\curvearrowright="3\msb@79
\mathchardef\digamma="0\msb@7A
\mathchardef\varkappa="0\msb@7B
\mathchardef\hslash="0\msb@7D
\mathchardef\hbar="0\msb@7E
\mathchardef\backepsilon="3\msb@7F
\def\Bbb{\ifmmode\let\next\Bbb@\else
 \def\next{\errmessage{Use \string\Bbb\space only in math mode}}\fi\next}
\def\Bbb@#1{{\Bbb@@{#1}}}
\def\Bbb@@#1{\fam\msbfam#1}

\mathchardef\cg="2\eufm@67
\mathchardef\cm="2\eufm@6D

\catcode`\@=12

\def\eps{{\epsilon}}

\def\<{\langle}
\def\>{\rangle}

\def\proof{\goodbreak\noindent{\bf Proof\quad}}

\def\text#1{{\rm #1}}
\def\note#1{}

\documentclass[11pt]{article}
\newtheorem{lemma}{Lemma}[section]
\newtheorem{propos}[lemma]{Proposition}
\newtheorem{example}[lemma]{Example}
\newtheorem{theorem}[lemma]{Theorem}
\newtheorem{cor}[lemma]{Corollary}

\newtheorem{defin}[lemma]{Definition}

\begin{document}
\baselineskip 22pt

\hfill Swan Preprint 02-1

\begin{center} {\LARGE POINTWISE BOUNDED ASYMPTOTIC 
MORPHISMS AND THOMSEN'S NON-STABLE k-THEORY}
\end{center}

\bigskip \begin{center} {{\Large E.J. Beggs} \\ \bigskip
Dept of Mathematics\\University of Wales Swansea\\
SA2 8PP U.K.\\ \medskip E.J.Beggs@swansea.ac.uk}
\end{center}

\subsection*{Abstract} In this paper I  show that pointwise bounded asymptotic morphisms
between separable metrisable locally convex *-algebras induce continuous maps between
the quasi-unitary groups of the algebras,
provided that the algebras support a certain amount of functional calculus.
 This links the asymptotic morphisms directly
to Thomsen's non-stable definition of $k$-theory in the $C^*$ algebra case.
A result on composition of asymptotic morphisms is also given.

\section{Introduction}
Thomsen defined a non-stable version of $K$-theory, called
 $k$-theory, in terms of the topology of the unitary group of a $C^*$ algebra \cite{Thom}. 
Here we link the pointwise bounded asymptotic morphisms (PBAMs) on separable
metrisable locally convex (SMLC) algebras to the topology of their unitary groups (for *-algebras). 
For other algebras, it may be more appropriate to consider the group of invertible elements,
for which a similar approach could be taken, but in this paper we concentrate on 
*-algebras and the unitary group. 

To avoid problems with algebras which do not contain a unit, Thomsen used the following 
procedure,
following Palmer \cite{Pal1,Pal2}.
On an algebra $A$ we define an associative binary operation $a\bullet b=a+b+a\,b$,
which has 2-sided unit $0\in A$. 
Note that if we had a unit $1$ for algebra multiplication,
 we would have $(1+a)(1+b)=1+a\bullet b$. By this means it can
be seen that this construction is really standard, though phrased in
slightly unusual terms. We define ${\rm gl}(A)$, the quasi-invertibles in 
an algebra $A$, and ${\cal U}(A)$, the quasi-unitaries, by
\[
{\rm gl}(A)\,=\, \{a\in A: \exists a'\in A\ \ a\bullet a'=a'\bullet a=0\}\ ,\ 
{\cal U}(A)\,=\, \{u\in A: u\bullet u^*=u^*\bullet u=0\}\ .
\]
It is easily seen that
both of these are groups under the $\bullet$ operation, with identity $0$. Thomsen's
non-stable $k$-groups for a $C^*$ algebra $E$ are defined as $k_n(E)=\pi_{n+1}({\cal U}(E))$.

The idea of asymptotic morphisms for $C^*$ algebras was introduced by Connes and Higson
 \cite{ConHig} for $E$-theory, and later considered in the non-stable $C^*$ case
by D\u ad\u arlat \cite{Dad}. In \cite{BegAsym} PBAMs were introduced to generalise the idea
of asymptotic morphism to SMLC algebras, and examples were given. The 
 definition of a PBAM will be given later in this paper. 
The general idea of all these constructions is to have a family of maps
$f_t:A\to B$ between algebras indexed by $t\in[0,\infty)$, where $f_t$ becomes 
more like an algebra map as $t\to\infty$. 

Under a certain assumption about functional calculus on SMLC *-algebras A and B, we show that
a PBAM $f:A\times [0,\infty)\to B$ induces a unique topological homotopy class ${\cal U}(f)\in
[{\cal U}(A),{\cal U}(B)]$. In the $C^*$ algebra case,
this gives a map $f_*:k_*(A)\to k_*(B)$ by composition,
$f_*[h]=[{\cal U}(f)\circ h]$. 
 The question naturally arises of whether ${\cal U}$ is a functor from
SMLC *-algebras and PBAMs to topological spaces. The answer is no, as SMLC *-algebras and PBAMs
themselves do not form a category. However there is a compatibility condition between
PBAMs $f:A\times [0,\infty)\to B$ and $g:B\times [0,\infty)\to C$ which ensures that
there is a composition $g\circ_\phi f:A\times [0,\infty)\to C$, defined up to homotopy of PBAMs. 
We shall see that this same condition also implies that ${\cal U}(g\circ_\phi f)=
{\cal U}(g)\circ{\cal U}(f)\in [{\cal U}(A),{\cal U}(C)]$. 
I choose to retain the compatibility condition rather than restrict to a category where 
composition is automatic (and there are such categories) in order not to unnecessarily
restrict the results.

\section{Asymptotic results and definitions on metric spaces}

\begin{lemma}\label{Y} Let $X$ be a separable metric space. 
Suppose that there is a function $\Psi:X\times X\times [0,\infty)\times [0,\infty)\to
\{{\rm true},{\rm false}\}$ with the following property: 
For all $x\in X$ there is a $\delta(x)>0$ and a $Q(x)\ge 0$ so that for all
$q\ge Q$ there is a $R(x,q)\ge 0$ and an $\epsilon(x,q)>0$ for which
\[
\forall y\in X\ \ \forall t\ge Q\ \ \forall s\ge R\quad
[\ t\in (q-\epsilon,q+\epsilon)
\ {\rm and}\ d(x,y)<\delta\ ]\ \Rightarrow\ \Psi(x,y,s,t)\ .
\]
Then there are continuous functions $\alpha:X\to[0,\infty)$ and $\phi:[0,\infty)
\to [0,\infty)$ (increasing) so that 
\[
\forall y\in X\ \ \exists x\in X\ \ \forall t\ge \alpha(y)\ \ 
\forall s\ge \phi(t)\quad [d(x,y)<\delta(x)\ {\rm and}\ \Psi(x,y,s,t)]\ .
\]
\end{lemma}
\proof {\bf Part 1:} Given $x\in X$ and $r\ge Q(x)$, there is an open cover
$(q-\epsilon(x,q),q+\epsilon(x,q))$ of the interval $[Q,r]$ indexed by
$q\in [Q,r]$. By compactness of $[Q,r]$ there is a finite subcover
indexed by $q_1,\dots,q_n$. 
Let $S(x,r)$ be the maximum of the $R(x,q_i)$. 
Then 
\[
\forall y\in X\ \ \forall t\in [Q,r]\ \ \forall s\ge S\quad
 d(x,y)<\delta\ \Rightarrow\ \Psi(x,y,s,t)\ .
\]
{\bf Part 2:} Take a cover of $X$ by open balls center $x\in X$ radius $\delta(x)$. 
There is a countable subcover of $X$ by open balls $B_i$ (integer $i\ge 0$) center $x_i$ radius
$\delta(x_i)$. 
Now define an increasing sequence $\phi_i\in[0,\infty)$ by using the 
recursive inequalities
\[
\forall j\in \{0,\dots i\}\quad{\rm if}\  i\ge Q(x_j)\ {\rm then}\ \phi_i\ge S(x_j,i+1)\ .
\]
Construct the graph of the continuous increasing function $\phi:[0,\infty)\to [0,\infty)$ by
joining the dots $(0,\phi_0)$, $(1,\phi_1)$, \dots with straight lines.
Take a continuous partition of unity $\theta_j:X\to[0,1]$
with locally finite supports, where $\theta_j$ has support in $B_j$, and define
\[
\alpha(x)\ =\ \sum_j (1+\max\{Q(x_j),j\})\, \theta_j(x)\ .
\]
{\bf Part 3:} 
For a given $y\in X$ we choose the minimum value of 
$\max\{Q(x_j),j\}$ over the finite number of indices $j$
 for which $\theta_j(y)\neq 0$. If this minimum
is $\max\{Q(x_k),k\}$ where $\theta_k(y)\neq 0$, it follows that $\alpha(y)\ge 1+
\max\{Q(x_k),k\}$ and $y\in B_k$. 
Now if $t\ge 1+\max\{Q(x_k),k\}$ then there is an $i\in\Bbb N$ with $i+1\ge t\ge i\ge Q(x_k)$
and $i\ge k$. 
Now if $s\ge \phi(t)$ then $s\ge \phi_i$, so $s\ge S(x_k,i+1)$. It follows that $\Psi(x_k,y,s,t)$. 
\quad$\square$

\begin{lemma}\label{YY} For $n\in\Bbb N$, let $X_n$ be a separable metric space. 
Suppose that there are functions $\Psi_n:X_n\times X_n\times [0,\infty)\times [0,\infty)\to
\{{\rm true},{\rm false}\}$ (for $n\in\Bbb N$) which satisfy the property
in the statement of lemma \ref{Y}, using functions
$\delta_n$, $Q_n$, $R_n$ and $\eps_n$. Then there is a continuous
increasing function $\psi:[0,\infty)
\to [0,\infty)$ and continuous functions $\beta_n:X_n\to[0,\infty)$ so that
\[
\forall n\in\Bbb N\ \ \forall y\in X_n\ \ \exists x\in X_n\ \ \forall t\ge \beta_n(y)\ \ 
\forall s\ge \psi(t)\quad [\ d(x,y)<\delta_n(x)\ {\rm and}\ \Psi_n(x,y,s,t)\ ]\ .
\]
\end{lemma}
\proof By using lemma \ref{Y} on each of the $\Psi_n$, we get continuous functions 
$\alpha_n:X_n\to[0,\infty)$ and $\phi_n:[0,\infty)
\to [0,\infty)$. Now we make a new continuous increasing function $\psi:[0,\infty)
\to [0,\infty)$ which satisfies the condition that if $t\ge n$ then $\psi(t)\ge \phi_n(t)$. 
Now set $\beta_n(x)=\max\{n,\alpha_n(x)\}$.\quad$\square$

\begin{defin}\label{sac} For metric spaces $X$ and $Y$, the function
$F:X\times[0,\infty)\to Y$ is said to be strongly asymptotically continuous
if, for all $\epsilon>0$ and all $x\in X$ there is an $\eta(x,\epsilon)>0$
and a $P(x,\epsilon)\ge 0$ so that
for all $x'\in X$ with $d(x,x')<\eta$ we have $d\big(F_t(x),F_t(x')\big)<\epsilon$
for all $t\ge P$.\end{defin}

\begin{lemma}\label{saccy} For metric spaces $X$ and $Y$, suppose that the function
$F:X\times[0,\infty)\to Y$ is continuous and strongly asymptotically continuous.
Then, given
 $\epsilon>0$ and $x\in X$, there are $\overline\eta(x,\epsilon)>0$ and 
 $\overline P(x,\epsilon)\ge 0$ so that, given $q\ge \overline P$
there is a $\overline \delta(x,\eps,q)>0$ so that
for all $x'\in X$ and all $t\ge 0$,
\[
[\ |t-q|<\overline\delta\ {\rm and}\ d(x,x')<\overline\eta\ ]\ \Rightarrow\ 
d\big(F_q(x),F_t(x')\big)<\epsilon\ .
\]
\end{lemma}\proof From definition \ref{sac} we take $\overline\eta=\eta(x,\epsilon/2)>0$
and $P(x,\epsilon/2)\ge 0$, and set $\overline P=P+1$. 
By continuity of $F$, given $q\ge \overline P$ there is a $1>\overline\delta(x,\eps,q)>0$ so that
$|t-q|<\overline \delta$ implies $d\big(F_q(x),F_t(x)\big)<\epsilon/2$. As $|t-q|<\overline \delta$
implies $t\ge P$, we have $d\big(F_t(x),F_t(x')\big)<\epsilon/2$ for all $x'\in X$ with
$d(x,x')<\overline\eta$. \quad$\square$

\section{Some functional calculus}

\begin{defin} \label{isrp}
Let $A$ be a topological *-algebra, and
take $A_{sad}$ to be the subset of self adjoint elements, i.e.\ those $a\in A$ for which
$a^*=a$. $A$ is said to have the inverse square root property if
there is a convex open set $V_A\subset A_{sad}$ containing $0\in A$ and a continuous function
$\vartheta:V_A\to A_{sad}$ so that

1)\quad $a\bullet\vartheta(a)=\vartheta(a)\bullet a$ for all $a\in V_A$,

2)\quad $a\bullet (\vartheta(a)\bullet \vartheta(a))=0$ for all $a\in V_A$,

3)\quad for all $a\in V_A$,  $a$ is $\bullet$-invertible,

4)\quad $\vartheta(0)=0$.
\end{defin}

\begin{lemma}\label{gp1} If $a\in A$ has $a^*\bullet a\in V_A$
and $a\bullet a^*\in V_A$, then $a\in{\rm gl}(A)$ and $a\bullet \vartheta(a^*\bullet a)
\in {\cal U}(A)$.
\end{lemma} 
\proof Set $u=a\bullet \vartheta(a^*\bullet a)$. Then
\[
u^*\bullet u\,=\, \vartheta(a^*\bullet a)\bullet (a^*\bullet a)\bullet \vartheta(a^*\bullet a)
\,=\, (a^*\bullet a)\bullet (\vartheta(a^*\bullet a)\bullet \vartheta(a^*\bullet a))\,=\,0\ .
\]
Since $a^*\bullet a$ and $a\bullet a^*$ are $\bullet$-invertible, there are $b\in A$ and $c\in A$
so that
\[
a^*\bullet a\bullet b\,=\,b\bullet a^*\bullet a\,=\,0\quad
{\rm and}\quad 
a\bullet a^*\bullet c\,=\,c\bullet a\bullet a^*\,=\,0\ .
\]
Now
\[
(a^*\bullet a)\bullet(b\bullet a^*)\,=\,
0\bullet a^*\,=\, a^*\,=\,a^*\bullet 0\,=\,a^*\bullet(a\bullet a^*\bullet c)
\,=\, (a^*\bullet a)\bullet(a^*\bullet c)\ ,
\]
so since $a^*\bullet a$ is $\bullet$-invertible we deduce that $b\bullet a^*=a^*\bullet c$
is the $\bullet$-inverse of $a$. The formula 
\[
((a^*\bullet a)\bullet \vartheta(a^*\bullet a))\bullet \vartheta(a^*\bullet a)\,=\, 0
\]
shows that $\vartheta(a^*\bullet a)$ is $\bullet$-invertible, so $u$ is $\bullet$-invertible, and so
$u\bullet u^*=0$, i.e.\ $u\in  {\cal U}(A)$.\quad$\square$

\begin{example} Any Banach *-algebra $A$ has the inverse square root property.
 This can be seen by
taking $V_A\,=\,\{a\in A_{sad}:|a|<\frac12\}$, and using functional calculus with the function
$\vartheta(x)=(1+x)^{-1/2}-1$. The easiest way to see that this works is to observe that
$\vartheta$ has a Taylor series with radius of convergence 1, and the fact that there is no 
constant term in the series means that we do not require $A$ to have a unit. 
\end{example}

\begin{example} The algebra $C^\infty(M,\Bbb C)$ of smooth complex valued
 functions on a compact manifold $M$ (with the 
smooth topology) has the inverse square root property. 
We set $V\,=\,\{f\in C^\infty(M,\Bbb R):\forall m\in M\ |f(m)|<\frac12\}$,
and $\vartheta(f)(m)=(1+f(m))^{-1/2}-1$.
\end{example}

\section{Pointwise bounded asymptotic morphisms}
Suppose that both $A$ and $B$ are separable
metrisable locally convex *-algebras (SMLC *-algebras for short).
Recall that a metrisable locally convex topological vector space has topology defined
by a countable number of seminorms, which we write $|.|_0$, $|.|_1$, etc. 
We can change any of the seminorms to a stronger continuous seminorm
 without altering the topology. It will be convenient to use this fact to assume that
the seminorms are non-decreasing, i.e.\ that $|b|_n\le|b|_{n+1}$ for all $n\in\Bbb N$ and $b\in B$.
  In a SMLC *-algebra $B$ we assume that 
the multiplication and star operations are continuous.
It will also be convenient to assume, by strengthening seminorms where appropriate, that
for all $b,b'\in B$ and $n\in\Bbb N$,
$|b^*|_n\le |b|_{n+1}$ and $|bb'|_n\le |b|_{n+1}|b'|_{n+1}$.
This implies the following result for $b,c\in B$, which is used later:
\begin{eqnarray}\label{useful45}
|b^*\bullet b-c^*\bullet c|_n\, \le\, 2\,|b-c|_{n+2}\ (1+|b|_{n+2}+|b-c|_{n+2})\ .
\end{eqnarray}

\begin{defin}\label{pb} The map $f:A\times [0,\infty)\to B$ is called pointwise asymptotically
bounded if, for all $a\in A$, the set $\{f_t(a):t\ge 0\}$ is bounded in $B$. 
In effect this means that for every $n\in\Bbb N$ the subset of the reals
$\{|f_t(a)|_n:t\ge 0\}$ is bounded above.
\end{defin}

\begin{defin}  If $A$ and $B$ are {\rm SMLC} $*-$algebras,
call a continuous function $f:A\times[0,\infty)\to B$ a pointwise bounded asymptotic
morphism (PBAM for short) if

a)\quad $f:A\times[0,\infty)\to B$ is pointwise asymptotically bounded
 (see \ref{pb}),

b)\quad $f$ is strongly asymptotically continuous (see \ref{sac}),

c)\quad $f_t(a)^*-f_t(a^*)\to 0$ as $t\to\infty$ for all $a\in A$,

d)\quad $\lambda\, f_t(a)-f_t(\lambda\, a)\to 0$ as $t\to\infty$ for all $a\in A$ and all
 $\lambda\in\Bbb C$,

e)\quad $f_t(a)+f_t(a')-f_t(a+a')\to 0$ as $t\to\infty$ for all $a,a'\in A$,

f)\quad $f_t(a)\,f_t(a')-f_t(a\,a')\to 0$ as $t\to\infty$ for all $a,a'\in A$.
\end{defin}

\medskip 
Note that the additive condition $(e)$ implies that $f_t(0)\to 0$ as $t\to\infty$.

\section{Pointwise bounded asymptotic morphisms and quasi-unitary groups}
In this section we suppose that $A$ and $B$ are SMLC *-algebras,
that $B$ has the inverse square root property (see \ref{isrp}), and that
$f:A\times[0,\infty)\to B$ is a 
PBAM. 

\begin{lemma}\label{bas}  For every $u\in {\cal U}(A)$ there
is a $T(u)\ge 0$ and an open neighbourhood $O(u)$ of $u$ in $A$ so that
\[
\forall t\ge T(u)\quad\forall v\in O(u)\quad 
f_t(v)\bullet f_t(v)^*\in V_B\ {\rm and}\ f_t(v)^*\bullet f_t(v)\in V_B\ .
\]
\end{lemma}
\proof First note that, for the given $u\in {\cal U}(A)$, 
$f_t(u)\bullet f_t(u)^*\to 0$ and $f_t(u)^*\bullet f_t(u)\to 0$
as $t\to\infty$. To see this, apply $f_t$ to $u^*\bullet u=0$
to get
$f_t(u)\,+\, f_t(u^*)\,+\, f_t(u^*)\,f_t(u)\to 0$  as $t\to\infty$. 
Now we use $f_t(u)^*-f_t(u^*)\to 0$  as $t\to\infty$ together with the fact that
$f_t(u)$ is bounded for $t\ge 0$ to show that $f_t(u)^*\bullet f_t(u)\to 0$. Similarly for 
$f_t(u)\bullet f_t(u)^*\to 0$. 

As $V_B$ is an open neighbourhood of $0$ in $B_{sad}$, there is
an $n\ge 0$ and an $\epsilon>0$ so that
$
\{ b\in B_{sad}:|b|_n<2\epsilon\}\subset V_B
$. 
As $f$ is pointwise bounded, there is a $K\ge 0$ so that $|f_t(u)|_{n+2}\le K$
for all $t$. Choose $\delta>0$ so that $2\delta(1+K+\delta)<\epsilon$. 
By strong asymptotic continuity there is a neighbourhood $O(u)$ of $u$
in $A$ and a $P\ge 0$ 
so that $|f_t(v)-f_t(u)|_{n+2}\le \delta$
for all $t\ge P(u)$ and all $v\in O(u)$. Then by (\ref{useful45}),
$|f_t(v)^*\bullet f_t(v)-f_t(u)^*\bullet f_t(u)|_n \,<\,\epsilon$. 
 Likewise we get 
$|f_t(v)\bullet f_t(v)^*-f_t(u)\bullet f_t(u)^*|_n<\,\epsilon$. Finally choose
$T\ge P$ so that $|f_t(u)^*\bullet f_t(u)|_n<\epsilon$ and $|f_t(u)\bullet f_t(u)^*|_n<\epsilon$
for all $t\ge T$. \quad $\square$

\begin{lemma}\label{alph} There is a continuous
function $\alpha:{\cal U}(A)\to [0,\infty)$ so that
\[
\forall v\in {\cal U}(A)\quad \forall t\ge \alpha(v)\quad
f_t(v)\bullet f_t(v)^*\in V_B\ {\rm and}\ f_t(v)^*\bullet f_t(v)\in V_B\ .
\]
\end{lemma}
\proof From lemma \ref{bas} take a countable subcover $O(u_n)$ ($n\in\Bbb N$) of ${\cal U}(A)$, and a 
locally finite continuous partition of unity $\theta_n$ on ${\cal U}(A)$, with 
the support of $\theta_n$ contained in $O(u_n)$. Now set $\alpha=\sum T(u_n)\,\theta_n$. 
\quad$\square$

\begin{propos} There is a well defined
homotopy class ${\cal U}(f)\in [{\cal U}(A),{\cal U}(B)]$ with representative function
$\tilde f_\alpha(v)\,=\, f_{\alpha(v)}(v)\bullet\vartheta(f_{\alpha(v)}(v)^*\bullet
f_{\alpha(v)}(v))$, where $\alpha$ is the
function in lemma \ref{alph}. Here well defined means
 independent of the choice of $\alpha$. 
\end{propos} 
\proof Note that the function $\tilde f_\alpha$ takes values in
${\cal U}(B)$ by lemma \ref{gp1}. 
Suppose that we have $\gamma:{\cal U}(A)\to[0,\infty)$ with 
$\gamma(u)\ge \alpha(u)$ for all $u\in {\cal U}(A)$. Then there is a homotopy
$H:{\cal U}(A)\times [0,1]\to {\cal U}(B)$ connecting $\tilde f_\alpha$ and $\tilde f_\gamma$
defined by $H(v,p)=\tilde f_{p\alpha+(1-p)\gamma}(v)$. 
Now if $\beta$ is another choice of function in lemma \ref{alph}, we see that
$\tilde f_\alpha$ and $\tilde f_\beta$ are homotopic, as they are both 
homotopic to $\tilde f_{{\rm max}\{\alpha,\beta\}}$.\quad$\square$

\section{Homotopy invariance}
In this section we suppose that $A$ and $B$ are SMLC *-algebras,
and that $f$ and 
$g:A\times[0,\infty)\to B$ are two PBAMs. 

\begin{defin} The PBAMs
 $f$ and $g$ are said to be pointwise bounded asymptotically
 homotopy equivalent (PBA homotopic for short) if there is a PBAM
 $h:A\times[0,\infty)\to C([0,1],B)$ so that 
$h_t(a)(0)=f_t(a)$ and $h_t(a)(1)=g_t(a)$.
(Use the supremum metric on $C([0,1],B)$, the continuous functions from
$[0,1]$ to $B$.)
\end{defin}

\begin{propos} If $f$ and $g$ are PBA
 homotopic then ${\cal U}(f)={\cal U}(g)\in [{\cal U}(A),{\cal U}(B)]$.
\end{propos} 
\proof Take the homotopy to be $h:A\times[0,\infty)\to C([0,1],B)$.
By lemma \ref{alph} there are functions 
$\alpha$, $\beta$ and $\gamma:{\cal U}(A)\to[0,\infty)$
corresponding to $f$, $g$ and $h$ respectively. Define $\eta:{\cal U}(A)\to[0,\infty)$
by $\eta(u)=\max\{\alpha(u),\beta(u),\gamma(u)\}$. Now $\tilde f_\alpha$
is homotopic to $\tilde f_\eta$, and $\tilde g_\beta$
is homotopic to $\tilde g_\eta$. Finally $\tilde f_\eta$ and $\tilde g_\eta$ are
homotopic by the homotopy $H:{\cal U}(A)\times [0,1]\to {\cal U}(B)$ defined by
$H(v,p)\,=\, \tilde h_\eta(v)(p)$, using the fact that
${\cal U}(C([0,1],B))=C([0,1],{\cal U}(B))$. \quad$\square$

\section{Composition of pointwise bounded asymptotic morphisms}
In this section we suppose that $A$, $B$ and $C$ are SMLC *-algebras, and that
$f:A\times[0,\infty)\to B$ and 
$g:B\times[0,\infty)\to C$ are PBAMs. We will show that under certain condtitions
it is possible to compose $f$ and $g$. To save a lot of time we will
use the results on composing strongly asymptotic morphisms in \cite{BegAsym}.

\begin{defin} For PBAMs $f$ and 
$g$, and a continuous increasing function $\phi:[0,\infty)\to [0,\infty)$
we define $g\circ_\phi f:A\times[0,\infty)\to C$ by 
 $(g\circ_\phi f)_t(a)=g_{\phi(t)}(f_t(a))$. 
 The function $\phi$
 is called a valid reparameterisation for $f$ and $g$
if 

a)\quad $g\circ_\phi f$ is a PBAM, and

b)\quad for every continuous increasing function
 $\theta:[0,\infty)\to [0,\infty)$ with $\theta(t)\ge \phi(t)$ $\forall
t\in[0,\infty)$, $g\circ_\theta f$ is also a PBAM, and  is
PBA homotopic to $g\circ_\phi f$, using
 $h:A\times [0,\infty)\to C([0,1],B)$ given by $h_t(a)(p)=g_{p\,\phi(t)+(1-p)\,\theta(t)}
(f_t(a))$. 

\end{defin}

\begin{defin}\label{cref}
 The PBAMs $f$ and 
$g$ are said to satisfy conditions {\it C1} and {\it C3} if:

{\it C1)}: For all $a\in A$ and all $\nu>0$ there is a $\xi(a,\nu)>0$
and a $Q'(a,\nu)\ge 0$ so that, for all $t\ge Q'$ there is a $S'(t,a,\nu)\ge 0$ so that
\[
[\ d(f_t(a),b)<\xi\ \text{and}\  s\ge S'\ ]\ \Rightarrow\ 
d\big(g_s(f_t(a)),g_s(b)\big)<\nu\ .
\]

{\it C3)}: For all $a\in A$ and all $n\in\Bbb N$ there is a $Q_n(a)\ge 0$ and an $M_n(a)\ge 0$
so that, for all $t\ge Q_n(a)$ there is a $S_n(a,t)\ge0$ so that for 
$s\ge S_n$ we have $ |g_s(f_t(a))|_n\le M_n$.
\end{defin}

\begin{lemma}\label{poi}  If the PBAMs $f$ and 
$g$ satisfy conditions {\it C1} and {\it C3}, then the
property $\Psi_n:A\times A\times[0,\infty)\times[0,\infty)\to\{{\rm true},{\rm false}\}$
defined by $\Psi_n(x,y,s,t)$ being the logical value of $|g_s(f_t(y))|_n\,\le\, M_n(x)+1$
satisfies the conditions of lemma \ref{Y}.
\end{lemma}
\proof Choose $\nu>0$ so that for the translation invariant metric $d$ on $C$, 
$d(c,c')<\nu$ implies $|c-c'|_n<1$. By applying {\it C3} to $x\in A$, we get 
$Q_n(x)\ge 0$ and an $M_n(x)\ge 0$
so that, for all $q\ge Q_n(x)$ there is a $S_n(x,q)\ge 0$ so that for
$s\ge S_n$ we have $ |g_s(f_q(x))|_n\le M_n$.
By applying {\it C1}, we get $\xi(x,\nu)>0$
and a $Q'(x,\nu)\ge 0$ so that, for all $q\ge Q'$ there is a $S'(q,x,\nu)\ge 0$ so that
\[
[\ d(f_q(x),b)<\xi\ \text{and}\  s\ge S'\ ]\ \Rightarrow\ 
|g_s(f_q(x))-g_s(b)|_n\,<\,1\ .
\]
By lemma \ref{saccy} there are $\overline\eta(x,\xi)>0$ and $\overline P(x,\xi)\ge 0$
so that, for $q\ge \overline P$ there is a $\overline\delta(x,\xi,q)>0$ so that
for all $y\in A$ and all $t\ge 0$,
\[
[\ d(x,y)<\overline\eta\ {\rm and}\ |t-q|<\overline\delta\ ]\ \Rightarrow\ d(f_q(x),f_t(y))<\xi\ .
\]

Now set $\delta(x)=\overline\eta(x,\nu)$ and 
$Q(x)=\max\{Q_n(x),Q'(x,\nu),\overline P(x,\xi)\}$. Given $q\ge Q$ set
$\eps(x,q)=\overline\delta(x,\xi,q)$ and $R(x,q)=\max\{S'(q,x,\nu),S_n(x,q)\}$. 
\quad$\square$

\begin{theorem}\label{reparam} If the PBAMs $f$ and 
$g$ satisfy conditions {\it C1} and {\it C3}, then there is a
continuous increasing function $\phi:[0,\infty)\to [0,\infty)$
which is a valid reparameterisation for $f$ and $g$.
\end{theorem}
\proof First note that  ({\it C1} and {\it C3}) implies {\it C2} from \cite{BegAsym}
for the *-operation, scalar multiplication, addition
and algebra multiplication. 
To see this we use the inequalities
\begin{eqnarray}\label{C2help}
|g_s(f_t(a))^*-g_s(b)^*|_n &\le& |g_s(f_t(a))-g_s(b)|_{n+1}\ , \cr
|\lambda\,g_s(f_t(a))-\mu\,g_s(b)|_n &\le& |\lambda-\mu|\,|g_s(f_t(a))|_{n}\,+\, \cr
&&
(|\lambda|+|\lambda-\mu|)\,|g_s(f_t(a))-g_s(b)|_{n}\ , \cr
|g_s(f_t(a))+g_s(f_t(a'))-g_s(b)-g_s(b')|_{n} &\le& |g_s(f_t(a))-g_s(b)|_{n}\,+\,
|g_s(f_t(a'))-g_s(b')|_{n}\ ,\cr
|g_s(f_t(a))\,g_s(f_t(a'))-g_s(b)\,g_s(b')|_{n} &\le& |g_s(f_t(a))|_{n+1}\,
|g_s(f_t(a'))-g_s(b')|_{n+1}\ +\cr
&&  |g_s(f_t(a))-g_s(b)|_{n+1}\ \times\cr
&&(|g_s(f_t(a'))|_{n+1}+|g_s(f_t(a'))-g_s(b')|_{n+1})\ .
\end{eqnarray}
From \cite{BegAsym} we now know that there is a valid reparameterisation
(for strongly asymptotic morphisms) $\psi:[0,\infty)\to[0,\infty)$ 
so that $g\circ_\psi f$ is a {\bf strongly} asymptotic morphism. From \ref{poi}
and \ref{YY} there are functions $\beta_n:A\to[0,\infty)$ ($n\in\Bbb N$)
 and $\theta:[0,\infty)\to [0,\infty)$ so that, given $n\in\Bbb N$ 
and $y\in A$, there is an $x\in A$ so that
\[
t\ge \beta_n(y) \ {\rm and}\ s\ge\theta(t)\ \Rightarrow\ |g_s(f_t(y))|_n\le M_n(x)+1\ .
\]
Let $\phi(y)=\max\{\psi(y),\theta(y)\}$. Then $g\circ_\phi f$ is a PBAM, and if $\chi:[0,\infty)
\to [0,\infty)$ has $\chi(t)\ge \phi(t)$ $\forall t\in [0,\infty)$, then 
the map $h:A\times [0,\infty)\to C([0,1],B)$ given by $h_t(a)(p)=g_{p\,\phi(t)+(1-p)\,\chi(t)}
(f_t(a))$ is a PBA homotopy between $g\circ_\phi f$ and $g\circ_\chi f$.\quad$\square$

\section{A conditional functoriality} 
Suppose that $A$, $B$ and $C$ are SMLC *-algebras, and that $B$ and $C$ have the 
 inverse square root property. We are given 
PBAMs $f:A\times[0,\infty)\to B$ and 
$g:B\times[0,\infty)\to C$ which satisfy {\it C1} and {\it C3}.

\begin{lemma}\label{ppp} Given $x\in {\cal U}(A)$, there is an $N(x)\ge 0$
and a $\nu(x)>0$ so that,
for $q\ge N$ there is an $L(x,q)\ge 0$ so that
\[
\forall z\in B\quad
[\ d(f_q(x),z)<\nu\ {\rm and}\ s\ge L\ ]\ \Rightarrow\ [\,g_s(z)^*\bullet g_s(z)\in V_C\
{\rm and}\ g_s(z)\bullet g_s(z)^*\in V_C\,]\ .
\]
\end{lemma}\proof {\bf Part 1:}
We shall only consider the $g_s(z)^*\bullet g_s(z)\in V_C$ part, the other is almost
identical.
Choose $\eps>0$ and $n\ge 0$ so that for $c\in C$, $|c|_n<\eps$ implies $c\in V_C$. 
Since $g_s(0)\to 0$ as $s\to\infty$, and using the
 strong asymptotic continuity of $g$, there is a $\eta>0$
and $P\ge 0$ so that
\begin{eqnarray}\label{ppp1}
[\ d(b,0)<\eta\ {\rm and}\ s\ge P\ ]\ \Rightarrow\ |g_s(b)|_n<\eps/4\ .
\end{eqnarray}
Given $x\in{\cal U}(A)$, by {\it C3} and (\ref{useful45})
 there is a $Q_{n+2}(x)\ge 0$ and and $M_{n+2}(x)\ge 0$ 
so that for any $q\ge Q_{n+2}$ there is a $S_{n+2}(x,q)$ so that for $s\ge S_{n+2}$,
\begin{eqnarray}\label{ppp2}
|g_s(z)^*\bullet g_s(z)-g_s(f_q(x))^*\bullet g_s(f_q(x))|_n &\le& 
2\,|g_s(z)-g_s(f_q(x))|_{n+2}\ \times  \cr
&& (\,1\,+\,M_{n+2}\,+\,|g_s(z)-g_s(f_q(x))|_{n+2})\ .
\end{eqnarray}
Now choose $\delta(x)>0$ so that $2\,\delta\,(1+M_{n+2}+\delta)<\epsilon/2$. 
By {\it C1} there is a $\xi(x)$ and a $Q'(x)$
so that, for all $q\ge Q'$ there is a $S'(x,q)$ so that
\begin{eqnarray}\label{ppp3}
[\ d(f_q(x),z)<\xi\ {\rm and}\ s\ge S'\ ]\ \Rightarrow\ |g_s(z)-g_s(f_q(x))|_{n+2}<\delta\ .
\end{eqnarray}
As $f$ is a PBAM we know that $f_t(x)^*\bullet f_t(x)-f_t(x^*\bullet x)\to 0$ as 
$t\to\infty$, and as $x\in{\cal U}(A)$ this means $f_t(x)^*\bullet f_t(x)\to 0$. 
Choose $N'(x)\ge 0$ so that $d(f_t(x)^*\bullet f_t(x),0)<\eta$ for all $t\ge N'$. 
Given $q\ge 0$, as $g$ is a PBAM we know that $g_s(f_q(x))^*\bullet g_s(f_q(x))-g_s(f_q(x)^*
\bullet f_q(x))\to 0$ as 
$s\to\infty$. Choose $S''(x,q)\ge 0$ so that
\begin{eqnarray}\label{ppp4}
s\ge S''\ \Rightarrow\ |\ g_s(f_q(x))^*\bullet g_s(f_q(x))-g_s(f_q(x)^*
\bullet f_q(x))   \,|_n<\eps/4\ .
\end{eqnarray}

\noindent
{\bf Part 2:} Set $N=\max\{Q_{n+2}(x),Q'(x),N'(x)\}$ and $\nu=\xi(x)$. Given
$q\ge N$ we set $L=\max\{P,S'(x,q),S''(x,q),S_{n+2}(x,q)\}$. Combining
(\ref{ppp2}) and (\ref{ppp3}), we have
\begin{eqnarray}\label{ppp5}
d(f_q(x),z)<\nu\ {\rm and}\ s\ge L\ \Rightarrow\ 
|g_s(z)^*\bullet g_s(z)-g_s(f_q(x))^*\bullet g_s(f_q(x))|_n<\eps/2\ .
\end{eqnarray}
As $q\ge N'$, we have $d(f_q(x)^*\bullet f_q(x),0)<\eta$, and 
putting this into (\ref{ppp1}) we get
\begin{eqnarray}\label{ppp6}
s\ge L\ \Rightarrow\ |g_s(f_q(x)^*\bullet f_q(x))|_n<\eps/4\ .
\end{eqnarray}
Combining this with (\ref{ppp4}) gives
\begin{eqnarray}\label{ppp7}
s\ge L\ \Rightarrow\ |g_s(f_q(x))^*\bullet g_s(f_q(x))|_n<\eps/2\ ,
\end{eqnarray}
which together with (\ref{ppp5}) gives the result.\quad$\square$

\begin{lemma}\label{qqq} 
For all $x\in{\cal U}(A)$, all $n\in\Bbb N$ and all $\eps>0$, there is
an $\eta_n(x,\eps)>0$ and a $T_n(x,\eps)\ge 0$ so that
for all $y\in A$ and all $p\in[0,1]$
\[
 t\ge T_n\ {\rm and}\ d(y,x)<\eta_n\ \Rightarrow\ |\vartheta_B(p(
f_t(y)^*\bullet f_t(y)))|_n<\eps\ .
\]
\end{lemma}
\proof The function $\vartheta_B:V_B\to B$ is continuous, so there is a
$m(\eps,n)\in\Bbb N$ and a $\delta(\eps,n)$ so that 
\[
\forall b\in A_{sad}\quad |b|_m<\delta\ \Rightarrow\ [\,b\in V_B
\ {\rm and}\ |\vartheta_B(b)|_n<\eps\,]\ .
\]
The problem then reduces to finding conditions under which
$|f_t(y)^*\bullet f_t(y)|_m<\delta$, and this is essentially contained in the proof of lemma 
\ref{bas}.\quad$\square$

\begin{defin} \label{defhat} For values $t\in[0,\infty)$, $u\in {\cal U}(A)$ and $p\in[0,1]$
for which $f_t(u)^*\bullet f_t(u)\in V_B$ and $f_t(u)\bullet f_t(u)^*\in V_B$, 
we define $
\hat f_{t,p}(u)\,=\, f_t(u)\bullet \vartheta(p(f_t(u)^*\bullet f_t(u)))
$.
\end{defin}

\begin{lemma}\label{bigp} The
 property $\Psi(x,y,s,t)$ defined as follows satisfies the conditions of lemma \ref{Y}. 
For $x,y\in X={\cal U}(A)$ and $s,t\ge 0$, $\Psi(x,y,s,t)$ is the logical value of the statement
\begin{eqnarray}\label{invbb}
f_t(y)^*\bullet f_t(y)\in V_B & {\rm and}& f_t(y)^*\bullet f_t(y)\in V_B
\ {\rm and}\cr \forall p\in[0,1]\quad[\  g_s(\hat f_{t,p}(y))^*\bullet g_s(\hat f_{t,p}(y))\in V_C
 & {\rm and}& g_s(\hat f_{t,p}(y))\bullet g_s(\hat f_{t,p}(y))^*\in V_C\ ]
\end{eqnarray}
\end{lemma}
\proof {\bf Part 1:}\quad By \ref{bas}, for every $x\in {\cal U}(A)$ there
is a $T(x)\ge 0$ and an open neighbourhood $O(x)$ of $x$ in $A$ so that
\[
\forall t\ge T(x)\quad\forall y\in O(x)\quad 
f_t(y)\bullet f_t(y)^*\in V_B\ {\rm and}\ f_t(y)^*\bullet f_t(y)\in V_B\ .
\]
By lemma \ref{ppp},
given $x\in {\cal U}(A)$, there is an $N(x)\ge 0$
and a $\nu(x)>0$ so that,
for $q\ge N$ there is an $L(x,q)\ge 0$ so that
\[
\forall z\in B\quad
[\ d(f_q(x),z)<\nu\ {\rm and}\ s\ge L\ ]\ \Rightarrow\ [\,g_s(z)^*\bullet g_s(z)\in V_C\
{\rm and}\ g_s(z)\bullet g_s(z)^*\in V_C\,]\ .
\]
We take $q\ge \max\{T,N\}+1$, $y\in O(x)\cap {\cal U}(A)$ and $|t-q|<1$ in what follows,
so $\hat f_{t,p}(y)$ (see \ref{defhat}) is defined for all $p\in[0,1]$. 

\noindent
{\bf Part 2:}\quad Take $n\in\Bbb N$ and $1>\chi>0$ so that in $B$, $|z-f_q(x)|_n<\chi$ implies
$d(f_q(x),z)<\nu$. Now consider the inequality
\begin{eqnarray}\label{ineq12}
|f_t(y)\bullet b-f_q(x)|_n\,\le\, |f_t(y)-f_q(x)|_{n+1}\,+\,
|b|_{n+1}\,(1+|f_q(x)|_{n+1}+|f_t(y)-f_q(x)|_{n+1})\ ,
\end{eqnarray}
where $b=\vartheta(p(f_t(y)^*\bullet f_t(y)))$, so that $f_t(y)\bullet b=\hat f_{t,p}(y)$. 
As $f$ is pointwise bounded there is an $M(x,n+1)$ so that $|f_q(x)|_{n+1}\le M$
for all $q\ge 0$. The conditions
\begin{eqnarray}\label{ineqs13}
|f_t(y)-f_q(x)|_{n+1}<\chi/2\quad {\rm and}\quad 
|\vartheta(p(f_t(y)^*\bullet f_t(y)))|_{n+1}<\chi/(4+2M)
\end{eqnarray}
 and $s\ge L(x,q)$ then imply $\Psi(x,y,s,t)$. 

\noindent
{\bf Part 3:}\quad From \ref{saccy} and \ref{qqq} there are $T_{n+1},\overline P\ge 0$ 
and $\eta_{n+1},\overline\eta>0$ (all depending on $x$) and, given $q\ge \overline P$,
a $\overline\delta>0$ (depending on $x$ and $q$) so that the conditions in (\ref{ineqs13})
are true if $t\ge T_{n+1}$, $|t-q|<\overline\delta$ and 
$d(x,y)<\min\{\eta_{n+1},\overline\eta\}$. Now define $Q(x)=\max\{T_{n+1},\overline P,T,N\}+1$,
and choose $\min\{\eta_{n+1},\overline\eta\}>\delta(x)>0$ so that $d(x,y)<\delta$
implies $y\in O(x)$. 
Given $q\ge Q$, set $R(x,q)=L(x,q)$ and $\eps(x,q)=\min\{\overline\delta,1\}$.
\quad$\square$

\begin{cor}\label{str}
There are
continuous functions $\gamma:{\cal U}\to[0,\infty)$ and $\theta:[0,\infty)\to [0,\infty)$
(increasing) so that
the formula
\[
r_{s,t,p}(u)\,=\, g_{s}(\hat f_{t,p}(u)) \bullet
\vartheta_C(\, g_{s}(\hat f_{t,p}(u))^*\bullet 
g_{s}(\hat f_{t,p}(u))\,)
\]
defines a continuous function 
\[
r:\{\ (u,s,t,p)\in {\cal U}(A)\times[0,\infty)\times[0,\infty)\times[0,1]:
t\ge\gamma(u)\ {\rm and}\ s\ge \theta(t)\ \}\to {\cal U}(C)\ .
\]
\end{cor}
\proof Use \ref{Y} on \ref{bigp}.\quad$\square$

\begin{theorem} ${\cal U}(g\circ f)={\cal U}(g)\circ {\cal U}(f)\in[{\cal U}(A),{\cal U}(C)]$. 
\end{theorem}
\proof By \ref{reparam}
 there is a valid reparameterisation $\phi:[0,\infty)\to [0,\infty)$ for $g$ and $f$.
Take the functions 
$\gamma:{\cal U}(A)\to[0,\infty)$ and $\theta:[0,\infty)\to [0,\infty)$ from lemma \ref{str},
and define $\psi:[0,\infty)\to [0,\infty)$ 
by $\psi(t)=\max\{t,\phi(t),\theta(t)\}$.
 As $f$, $g$ and $g\circ_\psi f$ are PBAMs, we can generate respective functions 
$\alpha:{\cal U}(A)\to[0,\infty)$, $\beta:{\cal U}(B)\to[0,\infty)$ and 
$\mu:{\cal U}(A)\to[0,\infty)$ from \ref{alph}. Define $\lambda:{\cal U}(A)\to[0,\infty)$
by $\lambda(u)=\max\{\alpha(u),\gamma(u)\}$
and $\omega:{\cal U}(A)\to[0,\infty)$ by
$\omega(u)=\max\{\alpha(u),\beta(\tilde f_{\lambda}(u)),\gamma(u),\mu(u)\}$.

Now the function $h:{\cal U}(A)\times[0,1]\to {\cal U}(C)$ defined by 
$
h_p(u)\,=\, r_{\psi(\omega(u)),\omega(u),p}(u)
$
gives a homotopy between $h_0(u)=(\widetilde{g\circ_\psi f})_\omega(u)$ and
\[
h_1(u)\,=\,  g_{\psi(\omega(u))}(\tilde f_\omega(u)) \bullet
\vartheta_C(\, g_{\psi(\omega(u))}(\tilde f_\omega(u))^*\bullet 
g_{\psi(\omega(u))}(\tilde f_\omega(u))\,)\ .
\]
By using $\tau,\zeta:[0,1]\times {\cal U}(A)\to[0,\infty)$ defined by
$\tau(p)(u)=(1-p)\,\omega(u)\,+\, p\,\lambda(u)$
and $\zeta(p)(u)=\max\{\psi(\tau(p)(u)),\, \beta(\tilde f_{\lambda}(u))\}$
 we define another homotopy 
$h':{\cal U}(A)\times[0,1]\to {\cal U}(C)$ by
\[
h'_p(u)\,=\,  g_{\zeta(p)(u)}(\tilde f_{\tau(p)}(u)) \bullet
\vartheta_C(\, g_{\zeta(p)(u)}(\tilde f_{\tau(p)}(u))^*\bullet 
g_{\zeta(p)(u)}(\tilde f_{\tau(p)}(u))\,)\ .
\]
This time $h'_0=h_1$ as $\psi(\omega(u))\ge\omega(u)\ge\beta(\tilde f_\lambda(u))$, and 
\[
h'_1(u)\,=\,  g_{\zeta(1)(u)}(\tilde f_{\lambda}(u)) \bullet
\vartheta_C(\, g_{\zeta(1)(u)}(\tilde f_{\lambda}(u))^*\bullet 
g_{\zeta(1)(u)}(\tilde f_{\lambda}(u))\,)\ .
\]
Define another homotopy $h'':{\cal U}(A)\times[0,1]\to {\cal U}(C)$ by
\[
h''_p(u)\,=\,  g_{\eta(p)(u)}(\tilde f_{\lambda}(u)) \bullet
\vartheta_C(\, g_{\eta(p)(u)}(\tilde f_{\lambda}(u))^*\bullet 
g_{\eta(p)(u)}(\tilde f_{\lambda}(u))\,)\ ,
\]
where $\eta(p)(u)=(1-p)\,\zeta(1)(u)\,+\,p\, \beta(\tilde f_{\lambda}(u))$. Now $h''_0=h'_1$ and
$h''_1=\tilde g_\beta\circ \tilde f_\lambda$.\quad$\square$

\end{document}